\documentclass{article}

\makeindex             

\begin{document}

\title{Characterizations of Probability Distributions via Bivariate
 Regression of Record Values}
\author{George P. Yanev \and
M. Ahsanullah \and M.I. Beg }

\maketitle

\begin{abstract}
Bairamov et al. (2005) characterize the exponential distribution
in terms of the regression of a function of a record value with
its adjacent record values as covariates. We extend these results
to the case of non-adjacent covariates. We also consider a more
general setting involving monotone transformations. As special
cases, we present characterizations involving weighted arithmetic,
geometric, and harmonic means.

\end{abstract}

 \section{Introduction and main results}

Let $X_1, X_2, \ldots$ be independent copies of a random variable
 $X$ whose distribution function is denoted by $F$. There is a number of studies on
 characterizations of F by means of regression relations of a
 function of one record value on one or two other record values. For a recent paper on the subject we refer to
 Pakes (2004) (see also Ahsanullah and Raqab (2006), Chapter 6). Denote upper record
 times by $L(1)=1$ and, for $n>1$,
 \[
 L(n)=\min \{j: \ j>L(n-1)\ \mbox {and}\ X_j>X_{L(n-1)}\},
 \]
 and the corresponding upper record value by $X(n)=X_{L(n)}$; see
 Nevzorov (2001). Gupta and Ahsanullah (2004) study the characterization of $F$ by means of the equation
 \[
 E[\psi(X(n))|X(n-k)=z]=\varphi(z),
 \]
 for $k=1$ and $k=2$, where the functions $\psi$ and $\varphi$
 satisfy certain regularity conditions.
 Bairamov et al. (2005) consider a
 characterization of the exponential distribution in terms of the regression on two adjacent
 record values
\[ E[\psi(X(n))|X(n-1)=u, X(n+1)=v ] \qquad (l_F<u<v<r_F),\]
where $l_F=\inf\{x:\ F(x)>0\}$ and $r_F=\sup\{x: \ F(x)<1\}$ are
 the left and right extremities of $F$, respectively.

 The aim of this paper is
 to extend the results given in Bairamov et al. (2005) by studying
 characterizations of $F$ in terms of the regression on two
 non-adjacent record values
\[ E[\psi(X(n))|X(n-k)=u, X(n+r)=v ] \qquad (l_F<u<v<r_F),\]
where $1\le k\le n-1$ and $r\ge 1$. Further on, for a given
function $h$, we adopt the notation \begin{equation}
 \label{M_not}
M(u,v)=\frac{h(v)-h(u)}{v-u}, \
_iM_j(u,v)=\frac{\partial^{i+j}}{\partial u^i\partial
v^j}\left(\frac{h(v)-h(u)}{v-u}\right) \quad (u\ne v),
\end{equation} as well as $_iM(u,v)$ and $M_j(u,v)$ for the $i$th
and $j$th partial derivative of $M(u,v)$ with respect to $u$ and
$v$, respectively. Bairamov et al. (2005) characterize the
exponential distribution as follows.

 {\bf Theorem (Bairamov et al. (2005)).}\ {\it Suppose $F$ is absolutely
continuous with density $f$, that $h$ is continuous in $[l_F,
r_F]$ and continuously differentiable in $(l_F, r_F)$, and that
almost everywhere in this open interval \begin{equation}
\label{ass} h'(x)\ne M(l_F, x). \end{equation} Then
\begin{equation} \label{cond1} E[h'(X(n))|X(n-1)=u,
X(n+1)=v]=M(u, v) \quad (l_F<u<v<r_F) \end{equation} holds if and
only if $l_F>-\infty$, $r_F=\infty$ and \[ F(x)=1-e^{-c(x-l_F)}
\quad (x\ge l_F), \] where $c>0$ is an arbitrary constant. }

Next result is a generalization of the above theorem to the case
of regression on a pair of non-adjacent record values. Namely,
(\ref{cond1}) is extended for $1\le k \le n-1$ and $r\ge 1$ to
\begin{eqnarray}
\label{eq2A} \lefteqn{\hspace{-1.2cm}E[h^{(k+r-1)}(X(n))|X(n-k)=u, X(n+r)=v]}\\
& & = \frac{(k+r-1)!}{(k-1)!(r-1)!}\ _{r-1}M_{k-1}(u, v) \qquad
(l_F<u<v<r_F). \nonumber
\end{eqnarray}
We consider two cases with respect to the spacings from the right
record value in the condition as follows: $X(n)$ is one spacing
away (there is a gap of size one); or $X(n)$ is two or more
spacings away (there is a gap of size two or more). The techniques
of the proofs in these two cases differ.

\newpage
{\bf Theorem 1.}\ {\it Suppose $F$ is absolutely continuous and
$h$ is continuous in $[l_F, r_F]$ and $h^{(k+r-1)}(x)$ is
continuous in $(l_F, r_F)$ for $1\le k \le n-1$ and  $r\ge 1$.

A. Let $r=1$, $1\le k \le n-1$ and \begin{equation} \label{ass2}
M_k(l_F,v)\ne 0 \qquad (l_F<v<r_F). \end{equation} Then
(\ref{eq2A}) holds if and only if \begin{equation}
\label{exp_type} F(x)=1-e^{-c(x-l_F)} \quad (x\ge l_F) \quad
\mbox{and}\quad l_F>-\infty, \ r_F=\infty,
\end{equation} where $c>0$ is an arbitrary constant.

B. Let $r\ge 2$, $2\le k\le n-1$ and $\ _rM_{k-1}(l_F,v)\ne 0$\
$(l_F<v<r_F)$. Then both (\ref{eq2A}) and
\begin{eqnarray} \label{eq2B} \lefteqn{\hspace{-0.9cm}E[h^{(k+r-1)}(X(n))|X(n-k+1)=u_2,
X(n+r)=v]} \\ & & = \frac{(k+r-2)!}{(k-2)!(r-1)!}\
_{r-1}M'_{k-2}(u_2, v) \qquad (l_F<u<u_2<v<r_F), \nonumber
\end{eqnarray}  where $M'(u,v)=[h'(v)-h'(u)]/(v-u)$, hold if and
only if (\ref{exp_type}) is true. }

{\bf Remarks.}\ (i) If $k=r=1$, then Theorem~1A coincides with
Theorem~1 in Bairamov et al. (2005) because, using
(\ref{main_first}) below, the assumption (\ref{ass}) can be
written as $M_1(l_F,v)\ne 0 $. (ii) The statement in Theorem 1A
holds true when $k=1$ and $r\ge 1$ as well (with $\ _rM(l_F,v)\ne
0$ instead of (\ref{ass2})) and can be proved along the same
lines, differentiating with respect to $u$ instead of $v$.

The following result is an extension of Theorem~2 in Bairamov et
al. (2005) to regression on a pair of non-adjacent covariates. As
we will see in the next section, it follows from Theorem~1
choosing $h(x)=x^{k+r}/(k+r)!$.

 {\bf Theorem 2.}\ A. Let $r=1$. For $1\le k\le n-1$ \begin{equation} \label{part_eq2A}
E[X(n)|X(n-k)=u, X(n+r)=v]=\frac{ru+kv}{k+1} \quad \quad
(l_F<u<v<r_F) \end{equation} if and only if the continuous r.v.
$X$ has the exponential distribution (\ref{exp_type}).

 B.  If $r\ge 2$ and $2\le k\le n-1$, then both (\ref{part_eq2A})
 and for $l_F<u<u_2<v<r_F$
 \[
 E[X(n)|X(n-k+1)=u_2, X(n+r)=v]=\frac{ru_2+(k-1)v}{r+k-1}\]
 hold if and only if the continuous r.v. $X$ has the
exponential distribution (\ref{exp_type}).

The proofs of Theorems 1 and 2 are given in Section 2. In
Section~3 we consider monotone transformations which extend the
results in the previous sections to a more general setting.
Illustrations are given in terms of characterizations involving
arithmetic, geometric, and harmonic means.

\section{Proofs of Theorems 1 and 2}

Further on we will need some recurrent relations for the
derivatives of $M(u,v)$ given in the following lemma.

{\bf Lemma 1.}\ Let $h(x)$ be a given function and for integer
$i,j\ge 1$ define $M(u,v)$, $_iM(u,v)$, $M_j(u,v)$, and
$_iM_j(u,v)$ as in (\ref{M_not}). If $h(x)$ has a continuous
derivative of order $\max\{i,j\}$ over the interval $(a,b)$, then
for $a<u<v<b$

\begin{equation} \label{main_first}
M_j(u,v)=\frac{h^{(j)}(v)-jM_{j-1}(u,v)}{v-u},\ \
_jM(u,v)=\frac{j\ _{j-1}M(u,v)-h^{(j)}(u)}{v-u} \end{equation} and
\begin{equation} \label{main} _iM_j(u,v)=\frac{i\
_{i-1}M_j(u,v)-j\ _iM_{j-1}(u,v)}{v-u}, \end{equation} where
 $M_1(u,v)$ and $_1M(u,v)$ are given in (\ref{main_first}) and
$_1M_1(u,v)=(M_1(u,v)\ -\ _1M(u,v))/(v-u)$.

 {\bf Proof.}\
It is not difficult to prove (\ref{main_first}) by induction. We
will proceed with the proof of (\ref{main}). One can check
(\ref{main}) for $i=1,2$ and $j=1,2$. Assume that (\ref{main})
holds for some $(i,j)$. Fixing $i$ we shall
prove it for $(i, j+1)$, i.e., \[ 
_{i-1}M_{j+1}(u,v)-(j+1)\ _iM_{j}(u,v)=(v-u) \ _iM_{j+1}(u,v) \]
Indeed, using the induction assumption, we have
\begin{eqnarray*} \lefteqn{
i\ _{i-1}M_{j+1}(u,v)-(j+1)\ _{i}M_{j}(u,v)}\\ &  = &
\frac{\partial}{\partial v}\left[ i\
_{i-1}M_{j}(u,v)-j\ _iM_{j-1}(u,v)\right]-\ _iM_j(u,v) \\
& = & \frac{\partial}{\partial v}\left[(v-u)\ _iM_j(u,v)\right]-\
_iM_j(u,v)\\
& = & (v-u)\ _{i}M_{j+1}(u,v).
 \end{eqnarray*}
 Similarly, assuming (\ref{main}) for an $i$ and fixed but arbitrary $j$,
 one can prove that it holds for $(i+1, j)$. The lemma is proved.

Denote $R(x)=-\ln(1-F(x))$. Using the Markov dependence of record
values, one can show (e.g., Ahsanullah (2004)) that the
conditional density of $X(n)$ given $X(n-k)=u$ $(1\leq k\le n-1)$
and $X(n+r)=v$ $(r\geq 1)$ is  \begin{equation} \label{cond_den}
\frac{\displaystyle (k+r-1)!}{\displaystyle
(k-1)!(r-1)!}\left[\frac{\displaystyle R(t)-R(u)}{\displaystyle
R(v)-R(u)}\right]^{k-1}\left[\frac{\displaystyle
R(v)-R(t)}{\displaystyle
R(v)-R(u)}\right]^{r-1}\frac{\displaystyle R'(t)}{\displaystyle
R(v)-R(u)},
\end{equation}
where $u<t<v$.

We will need the following lemma, which is of independent interest
as well.

{\bf Lemma 2.}\ {\it Let $h(x)$ be a continuous in $[l_F, r_F]$
function such that $h^{(k+r-1)}(x)$ is continuous in $(l_F, r_F)$
for $1\le k \le n-1$ and  $r\ge 1$.  If
\[ F(x)=1-e^{-c(x-l_F)} \quad (x\ge l_F) \quad \mbox{and}\quad
l_F>-\infty, \ r_F=\infty, \] where $c>0$ is an arbitrary
constant, then
\begin{eqnarray}
\label{eq2AL} \lefteqn{\hspace{-1.2cm}E[h^{(k+r-1)}(X(n))|X(n-k)=u, X(n+r)=v]}\\
& & = \frac{(k+r-1)!}{(k-1)!(r-1)!}\ _{r-1}M_{k-1}(u, v) \qquad
(l_F<u<v<r_F). \nonumber
\end{eqnarray}
}

 {\bf Proof.}\ It is not difficult to verify (\ref{eq2AL}) for
 $k=r=1$. Let us prove it for $r=1$ and any $1\le k\le n-1$, i.e.,
\begin{equation} \label{k=1}
 E[h^{(k)}(X(n))|X(n-k)=u, X(n+1)=v]=
kM_{k-1}(u, v) \end{equation} for $1\le k\le n-1$.  Assuming that
(\ref{k=1}) is true for $k=i$ \ $(1\le i\le n-1)$, we will prove
it for $k=i+1$. Indeed, making use of (\ref{cond_den}) with
$R(x)=c(x-l_F)$ and the induction assumption, we obtain
\begin{eqnarray*} \lefteqn{E[h^{(i+1)}(X(n))|X(n-i-1)=u,
X(n+1)=v] }\\ & = &
\frac{i+1}{(v-u)^{i+1}}\int_u^v h^{(i+1)}(t)(t-u)^{i}dt\\
& = & \frac{i+1}{(v-u)^{i+1}} \left[ h^{(i)}(v)(v-u)^i -i\int_u^v
h^{(i)}(t)(t-u)^{i-1}dt\right] \\
& = & \frac{i+1}{(v-u)^{i+1}} \left[ h^{(i)}(v)(v-u)^i -
(v-u)^iE[h^{(i)}(X(n))|X(n-i)=u,
X(n+1)=v]\right] \\
 & = &
\frac{i+1}{(v-u)^{i+1}} \left[ h^{(i)}(v)(v-u)^i
- i(v-u)^iM_{i-1}(u,v)\right]\\
& = & \frac{i+1}{v-u}\left[h^{(i)}(v)-iM_{i-1}(u,v)\right]\\
& = & (i+1)M_i(u,v),
 \end{eqnarray*}
where the last equality follows from (\ref{main_first}). This
proves (\ref{k=1}) for $k=i+1$ and thus (\ref{eq2AL}) is true  for
$r=1$ and any $1\le k\le n-1$. Similarly one can prove
(\ref{eq2AL}) for $k=1$ and any $r\ge 1$, i.e.,
\begin{equation}\label{r=1}
E[h^{(r)}(X(n))|X(n-1)=u, X(n+r)=v]=r\ _{r-1}M(u,v)
\end{equation}

To complete the proof of the lemma we need to show (\ref{eq2AL})
for $r\ge 2$ and $2\le k\le n-1$. Let us assume (\ref{eq2AL}) for
$r=j$ and any $2\le k\le n-1$. We will prove it for $r=j+1$ and
any $2\le k\le n-1$. Since the left-hand side of (\ref{eq2AL}) for
$r=j+1$ is
\begin{eqnarray*}
\lefteqn{E[h^{(k+j)}(X(n))|X(n-k)=u, X(n+j+1)=v]} \\
&= &\frac{(k+j)!(v-u)^{-(k+j)}}{(k-1)!j!}\int_u^v
h^{(k+j)}(t)(t-u)^{k-1}(v-t)^rdt,
\end{eqnarray*}
to prove (\ref{eq2AL}) for $r=j+1$ we need to show that for $2\le
k\le n-1$
\begin{eqnarray} \label{integral=M} I(k,j+1)& = & \int_u^v
h^{(k+j)}(t)(t-u)^{k-1}(v-t)^jdt \\
 & =  & (v-u)^{k+j}\ _jM_{k-1}(u,v) \nonumber
\end{eqnarray}
under the induction assumption that for any $2\le k\le n-1$
\begin{eqnarray} \label{integral=M2} I(k,j) & = & \int_u^v
h^{(k+j-1)}(t)(t-u)^{k-1}(v-t)^{j-1}dt \\
    & = & (v-u)^{k+j-1}\
_{j-1}M_{k-1}(u,v) \nonumber
\end{eqnarray}
Integrating by parts, we have for $2\le k\le n-1$
\begin{eqnarray*}
 I(k,j)
 & = & \int_u^vh^{(k+j-1)}(t)(t-u)^{k-1}(v-t)^{j-1}dt  \\
    & = & \frac{1}{j}\int_u^v h^{(k+j)}(t)(t-u)^{k-1}(v-t)^jdt
    \\
    & & +
    \frac{k-1}{j}\int_u^v h^{(k+j-1)}(t)(t-u)^{k-1}(v-t)^jdt
    \end{eqnarray*}
    \noindent Let $(x)_n$ be the falling factorial, i.e.,
$(x)_n=x(x-1)\ldots(x-n+1)$ $(n\ge 1)$ and $(x)_0=1$.
    After iterating, we have for $2\le k\le n-1$
\begin{eqnarray} \label{long_expr} \lefteqn{I(k,j+1)}\\
& = & jI(k,j)-(k-1)I(k-1,j+1) \nonumber \\
  & = & j\sum_{i=0}^{k-2}
(-1)^{i}(k-1)_{(i)}I(k-i,j) +
(-1)^{k-1}(k-1)_{(k-1)}\int_u^vh^{(j+1)}(t)(v-t)^{j}dt \nonumber
\end{eqnarray}
Observe that (\ref{r=1}) and (\ref{main_first}) lead to
\begin{eqnarray} \label{last_int}
\lefteqn{\int_u^vh^{(j+1)}(t)(v-t)^jdt}\\
 & = &
-h^{(j)}(u)(v-u)^j+j\int_u^vh^{(j)}(t)(v-t)^{j-1}dt \nonumber \\
    & = &
(v-u)^j\left\{-h^{(j)}(u)+E[h^{(j)}(X(n))|X(n-1)=u, X(n+j)=v]\right\}\nonumber \\
    & = &
    (v-u)^j[-h^{(j)}(u)+j\ _{j-1}M(u,v)] \nonumber \\
        & = &
 (v-u)^j[-h^{(j)}(u)+(v-u)\ _jM(u,v)+h^{(j)}(u)] \nonumber
 \\
    & = & (v-u)^{j+1}\ _jM(u,v) \nonumber
\end{eqnarray}
Using the induction assumption (\ref{integral=M2}) and
(\ref{last_int}) we write (\ref{long_expr}) as
\begin{eqnarray} \label{long_expr2}
\frac{I(k,j+1)}{(v-u)^{j+1}} & = &
j\sum_{i=0}^{k-2}(-1)^{i}(k-1)_{(i)}(v-u)^{k-2-i}\
_{j-1}M_{k-1-i}(u,v)\\
&  & -(-1)^{k-2}(k-1)_{(k-1)}\ _jM(u,v)\nonumber
\end{eqnarray}
Now, applying (\ref{main}) and iterating, we obtain
\begin{eqnarray*}
\lefteqn{\frac{I(k,j+1)}{(v-u)^{j+1}}}\\
 & = &
j\sum_{i=0}^{k-3}(-1)^{i}(k-1)_{(i)}(v-u)^{k-2-i}\
_{j-1}M_{k-1-i}(u,v)\\
    & & \hspace{1cm} +(-1)^{k-2}(k-1)_{(k-2)}\left[j\ _{j-1}M_1(u,v)-\
    _jM(u,v)\right] \\
    & = &
j\sum_{i=0}^{k-3}(-1)^{i}(k-1)_{(i)}(v-u)^{k-2-i}\
_{j-1}M_{k-1-i}(u,v) \\
& & \hspace{1cm} -(-1)^{k-3}(k-1)_{(k-2)}(v-u)\ _jM_1(u,v)\\
    & = &
j\sum_{i=0}^{k-4}(-1)^{i}(k-1)_{(i)}(v-u)^{k-2-i}\
_{j-1}M_{k-1-i}(u,v)\\
& & \hspace{1cm} +(-1)^{k-3}(k-1)_{(k-3)}(v-u)\left[j\
_{j-1}M_2(u,v)-\
    2\ _jM_1(u,v)\right]\\
 & = &
j\sum_{i=0}^{k-4}(-1)^{i}(k-1)_{(i)}(v-u)^{k-2-i}\
_{j-1}M_{k-1-i}(u,v)\\
& & \hspace{1cm} -(-1)^{k-4}(k-1)_{(k-3)}(v-u)^2\  _{j-1}M_2(u,v)\\
& & \cdots \\
& = & j(v-u)^{k-2}\ _{j-1}M_{k-1}-(k-1)(v-u)^{k-2}\ _{j-1}M_{k-2}(u,v)\\
& = & (v-u)^{k-1}\ _{j-1}M_{k-1}(u,v).
\end{eqnarray*}
This implies (\ref{integral=M}). Similarly assuming (\ref{eq2AL})
for $k=i$ $(2\le i\le n-2)$ and any $r\ge 2$ one can prove it for
$k=i+1$ and $r\ge 2$. The lemma is proved.

\subsection{Proof of Theorem 1A}

Assume (\ref{eq2A}). Setting $r=1$ in (\ref{cond_den}), we obtain
from (\ref{eq2A})
\[
M_{k-1}(u,v)[R(v)-R(u)]^k =\int_{u}^v
h^{(k)}(t)[R(t)-R(u)]^{k-1}R'(t)dt.
\]
Letting $u\to l_F^+$ and noting that the integrand is continuous
and $\lim_{u\to l_F^+}R(u)=\lim_{u\to l_F^+}(-\ln (1-F(u)))=0$ we
simplify to
\[
M_{k-1}(l_F,v)[R(v)]^k =\int_{l_F}^v
h^{(k)}(t)[R(t)]^{k-1}R'(t)dt.
\]
 Differentiating both sides of the above equation with
respect to $v$, we obtain
\[ kM_{k-1}(l_F, v)[R(v)]^{k-1}R'(v)+M_k(l_F,v)[R(v)]^k=
h^{(k)}(v)[R(v)]^{k-1}R'(v)
\]
Rearranging and taking into account (\ref{main_first}),
 \begin{eqnarray*} \frac{R'(v)}{R(v)} & = &
\frac{M_k(l_F,v)}{h^{(k)}(v)-kM_{k-1}(l_F,v)}\\
    & = &
\frac{M_k(l_F,v)}{(v-l_F)M_k(l_F, v)}\\
    & = &
    \frac{1}{v-l_F},
    \end{eqnarray*}
(provided that $M_k(l_F,v)\ne 0$) and hence (\ref{exp_type})
holds. It follows that $l_F>-\infty$ and $c>0$, and the continuity
of $F$ implies that $r_F=\infty$.

The converse statement follows from Lemma 2. The proof is
complete.

\subsection{Proof of Theorem 1B}\

Assume both (\ref{eq2A}) and (\ref{eq2B}) are true. Formula
(\ref{cond_den}) together with (\ref{eq2A}) imply
\begin{eqnarray*}
\lefteqn{_{r-1}M_{k-1}(u,v)\left[R(v)-R(u)\right]^{k+r-1}}\\
& = & \int_{u}^vh^{(k+r-1)}(t)\left[R(v)-R(t)\right]^{r-1}
\left[R(t)-R(u)\right]^{k-1}R'(t)dt.\end{eqnarray*} Since the
integrand is continuous, differentiating both sides of the above
equation with respect to $u$, we obtain
\begin{eqnarray} \label{diff1}
\lefteqn{_rM_{k-1}(u,v)\left[R(v)-R(u)\right]^{k+r-1}}\\
& & -(k+r-1)\left[R(v)-R(u)\right]^{k+r-2} R'(u)\
_{r-1}M_{k-1}(u,v) \nonumber \\
& &
=-(k-1)R'(u)\int_{u}^vh^{(k+r-1)}(t)\left[R(v)-R(t)\right]^{r-1}
\left[R(t)-R(u)\right]^{k-2}R'(t)dt. \nonumber
\end{eqnarray}
On the other hand, (\ref{eq2B}) and (\ref{cond_den}) lead to
\begin{eqnarray} \label{u2_cond}
\lefteqn{_{r-1}M'_{k-2}(u_2,v)\left[R(v)-R(u_2)\right]^{k+r-2}}\\
 = & & \int_{u_2}^vh^{(k+r-1)}(t)\left[R(v)-R(t)\right]^{r-1}
\left[R(t)-R(u_2)\right]^{k-2}R'(t)dt. \nonumber \end{eqnarray}
 Therefore, letting $u_2\to u^+$ in (\ref{u2_cond}) and rearranging terms, we write
 (\ref{diff1}) as
 \begin{equation} \label{eqn} \frac{\displaystyle R'(u)}{\displaystyle R(u)-R(v)} =
\frac{\displaystyle _rM_{k-1}(u,v)}{\displaystyle (k-1)\
_{r-1}M'_{k-2}(u,v)-(k+r-1)\ _{r-1}M_{k-1}(u,v)} \end{equation}
provided that the denominator in the right-hand side is not 0.
(This is equivalent to $\ _rM_{k-1}(u,v)\ne 0$, as we will see
below.) Since
 \begin{eqnarray*}
_{r-1}M'_{k-2}(u,v) & = & \ \frac{\partial^{k+r-3}}{\partial
u^{r-1}\partial v^{k-2}}\left[\frac{h'(v)-h'(u)}{v-u}\right]\\
    & = &
 \frac{\partial^{k+r-3}}{\partial
u^{r-1}\partial v^{k-2}} \left[M_1(u,v)+\ _1M(u,v)\right]\\
    & = &
\ _{r-1}M_{k-1}(u,v)+\ _{r}M_{k-2}(u,v), \end{eqnarray*} for the
denominator in (\ref{eqn}) we have \begin{eqnarray}
\lefteqn{\hspace{-0.5cm}(k-1)\ _{r-1}M'_{k-2}(u,v)-(k+r-1)\
_{r-1}M_{k-1}(u,v)}\label{denominator}\\
 & = &
(k-1)[\ _{r-1}M_{k-1}(u,v)+\ _{r}M_{k-2}(u,v)]-(k+r-1)\
_{r-1}M_{k-1}(u,v)
 \nonumber \\
    & = &
    (k-1)\ _{r}M_{k-2}(u,v)-r\
    _{r-1}M_{k-1}(u,v). \nonumber
\end{eqnarray} Finally, from (\ref{eqn})-(\ref{denominator}) and applying
(\ref{main}), we obtain \begin{eqnarray*} \frac{\displaystyle
R'(u)}{\displaystyle R(u)-R(v)} & = & \frac{\displaystyle \
_{r}M_{k-1}(u,v)}{ \displaystyle (k-1)\
    _{r}M_{k-2}(u,v) - r\ _{r-1}M_{k-1}(u,v) }\\
    & = & \frac{\displaystyle \
_{r}M_{k-1}(u,v)}
{\ _{r}M_{k-1}(u,v)(u-v)}\\
    & = &
    \frac{\displaystyle 1}{\displaystyle u-v}.
    \end{eqnarray*}
Integrating both sides with respect to $u$ from $l_F$ to $v$, we
obtain
\[
\ln[R(v)-R(l_F)]=\ln (v-l_F) + \ln c    \qquad (c>0)
\]
and thus $ R(v)=c(v-l_F)$ and (\ref{exp_type}) follows.

The converse statement in the theorem follows from Lemma 2.

\subsection{Proof of Theorem 2}

  Let $h(x)=x^{k+r}/(k+r)!$ and
  thus $h^{(k+r-1)}(x)=x$. We shall prove that, with this choice of $h$, (\ref{eq2A})
  becomes
\begin{equation} \label{ex1} E[X(n)|X(n-k)=u, X(n+r)=v]=\frac{ru+kv}{k+r} \quad
(l_F<u<v<r_F). \end{equation}
 Indeed,
\begin{eqnarray*}
M(u,v) & = & \frac{1}{(k+r)!}\frac{v^{k+r}-u^{k+r}}{v-u} \\
& = & \frac{v^{k+r-1}+\ldots +v^{k}u^{r-1}+v^{k-1}u^r+\ldots
+u^{k+r-1}}{(k+r)!}
\end{eqnarray*}
and differentiating $r-1$ and $k-1$ times with respect to $u$ and
$v$, we obtain \begin{eqnarray} \label{diff3A}
\lefteqn{\hspace{-2.5cm} \frac{(k+r-1)!}{(k-1)!(r-1)!}\
_{r-1}M_{k-1}(u,v)} \\ & &
=\frac{(k+r-1)!}{(r-1)!(k-1)!}\frac{(r-1)!k!v+r!(k-1)!u}{(k+r)!}\nonumber
\\
& & =\frac{ru+kv}{k+r}, \nonumber
\end{eqnarray} which proves (\ref{ex1}). Now, if $r=1$ (note that
$M_k(l_F,v)=l_F/(k+1)\ne 0$) the claim in Theorem 2A follows from
Theorem 1A. Similarly to (\ref{diff3A}) one can see that
(\ref{eq2B}) becomes
\begin{eqnarray} \label{diff3B}
 \lefteqn{\hspace{-2.5cm} E[X(n)|X(n-k+1)=u_2,
X(n+r)=v]}\\
& & =\frac{(k+r-2)!}{(k-2)!(r-1)!}\ _{r-1}M'_{k-2}(u_2,v)\nonumber \\
& & = \frac{ru_2+(k-1)v}{k+r-1}. \nonumber
 \end{eqnarray}
 Theorem 1B, (\ref{diff3A}) and (\ref{diff3B})
 imply Theorem 2B. The proof is complete.

\section{Monotone transformations and some particular cases}

In this section, following Bairamov et al. (2005), we give a
formal generalization of Theorem~1, involving a monotone
transformation of $X$. Let $Y$ be a random variable with
distribution function $G$. The corresponding upper record values
are denoted by $Y(n)$. The following extension of Theorem~1A
holds. The proof is similar to that of Theorem~3 in Bairamov et
al. (2005) and it is omitted here. Denote for $i,j\ge 0$
\[
_iM_j(T(s),T(t))=\frac{\partial^{i+j}}{\partial x^i\partial
y^j}\left(\frac{h(y)-h(x)}{y-x}\right)|_{x=T(s),\ y=T(t)} \qquad
(y\ne x).
\]
 {\bf Theorem 3.}\ Suppose that:

 (i) $Y$ has a continuous distribution function $G$ on
 $[l_G,r_G]$;

 (ii) the function $T$ is continuous and strictly increasing in
 $(l_G,r_G)$, $\tau=T(l_{G+})>-\infty$ and $T(r_G)=\infty$; and

 (iii) $h^{(k+r-1)}(x)$ is continuous in $(\tau, \infty)$ for $1\le k\le
 n-1$ and $r\ge 1$.

A. Let $r=1$ and $M_k(\tau,T(t))\ne 0 \ (l_G<t<r_G)$. Then for
$1\le k \le n-1$ and $l_G<s<t<r_G$
\begin{eqnarray} \label{eq3A}
E[h^{(k)}(T(Y(n)))|Y(n-k)=s, Y(n+1)=t]  & = & kM_{k-1}(T(s), T(t))
 \end{eqnarray} if and only if
\begin{equation} \label{gen_exp} G(y)=1-e^{-c[T(y)-\tau]} \quad (l_G<y<r_G) \end{equation}
where $c>0$ is an arbitrary constant.

B. Let $r\ge 2$ and $2\le k\le n-1$, and $\ _rM_{k-1}(\tau,
T(t))\ne 0$ where  $l_G<t<r_G$. Then both (\ref{eq3A}) and for
$l_G<s<s_2<t<r_G$
\begin{eqnarray*} \lefteqn{E[h^{(k+r-1)}(T(Y(n)))|Y(n-k+1)=s_2, Y(n+r)=t]} \\
& & = \frac{(k+r-2)!}{(k-2)!(r-1)!}\ _{r-1}M'_{k-2}(T(s_2), T(t))
\end{eqnarray*}
 hold if and only if (\ref{gen_exp}) is true.

{\bf Remark.}\ An analog of Theorem~3 when $T$ is a strictly
decreasing function holds as well; for the case $k=r=1$ see
Bairamov et al. (2005), Theorem~3.

Different choices of functions $h$ and $T$ in the above theorem
yield many characterization results. The corollary below gives a
characterization that involves a weighted arithmetic mean.

{\bf Corollary 1.}\ Let (i) and (ii) of Theorem~3 hold. If $1\le
k\le n-1$, then for a strictly increasing function $g$
\[
E[g(Y(n))|Y(n-k)=s, Y(n+1)=t]=\frac{kg(t)+g(s)}{k+1} \qquad
(l_G<s<t<r_G)
\]
holds if and only if \[ G(y)=1-e^{\displaystyle -c[g(y)-g(\tau)]}
\qquad (l_G<y<r_G),
\] where $c>0$ is an arbitrary constant.

{\bf Proof.}\ The corollary follows from Theorem 2A and Theorem~3
setting $T(y)=g(y)$ and $h(x) =x^{k+1}/(k+1)!$.

Next corollary presents characterizations involving a geometric
mean as a special case.

{\bf Corollary 2.}\ Let (i) and (ii) of Theorem~3 hold. If $1\le
k\le n-1$, then for a strictly decreasing function $g$ and
$l_G<s<t<r_G$ \begin{equation} \label{gen_geo} E[g(Y(n))|Y(n-k)=s,
Y(n+1)=t]=[g(t)]^{k/(k+1)}[g(s)]^{r/(k+1)}
\end{equation} holds if and only if
\[ G(y)=1-e^{\displaystyle -c\left\{ [g(y)]^{-1/(k+1)}-[g(\tau)]^{-1/(k+1)}\right\}} \qquad
(l_G<y<r_G),
\] where $c>0$ is an arbitrary constant.

 {\bf Proof.}\ We will show that if $h(x)=-x^{-1}$,
 then for $j=1,2,\ldots$
\begin{equation} \label{bigM} M_j(x,y)=(-1)^j\frac{j!}{xy^{j+1}}. \end{equation} Indeed,
one can check that $M_1(x,y)=-1/(xy^2)$. Assuming that
(\ref{bigM}) is true for $j$, we will prove it for $j+1$. Using
(\ref{main_first}) we have \begin{eqnarray*} M_{j+1}(x,y) & = &
\frac{1}{y-x}\left[
(-1)^{j+2}\frac{(j+1)!}{y^{j+2}}-(j+1)(-1)^j\frac{j!}{xy^{j+1}}\right]\\
& = & \frac{1}{y-x}\left[ (-1)^{j+1}\frac{(j+1)!(y-x)}{xy^{j+2}}\right]\\
& = & (-1)^{j+1}\frac{(j+1)!}{xy^{j+2}} \end{eqnarray*} and thus
(\ref{bigM}) follows by induction. Now, let us set for $1\le k\le
n-1$
\[
h(x)=\frac{(-1)^k}{k!x} \qquad \mbox{and thus} \qquad
h^{(k)}(x)=\frac{1}{x^{k+1}}.
\]

It is not difficult to see that, with this choice of $h$,
(\ref{eq3A}) and (\ref{bigM}) yield
\[
E\left[\frac{1}{[T(Y(n))]^{k+1}}|Y(n-k)=s, Y(n+1)=t\right]=
\frac{1}{T(s)[T(t)]^k}.
 \]
Setting  in the last equation $T(x)=[g(x)]^{-1/(k+1)}$, leads to
the statement of the corollary.

It is worth noting that if $k=r=1$, then the right-hand side in
(\ref{gen_geo}) is the geometric mean of $g(s)$ and $g(t)$.

Next corollary is a characterization in terms of a harmonic mean.

{\bf Corollary 3.}\ Let (i) and (ii) of Theorem~3 hold. For a
strictly decreasing function $g$,
\[
E[g(Y(n))|Y(n-1)=s, Y(n+1)=t]=\frac{2g(s)g(t)}{g(s)+g(t)} \qquad
(l_G<s<t<r_G).
\]
holds if and only if
\[ G(y)=1-e^{\displaystyle -c\left\{ [g(y)]^{-2}-[g(\tau)]^{-2}\right\}} \qquad
(l_G<y<r_G)
\] where $c>0$ is an arbitrary constant.

{\bf Proof.}\ The result follows from Theorem~3 setting
$h(x)=2x^{1/2}$ and $T(y)=[g(y)]^{-2}$.

Finally, let us note that Theorem~3 and its corollaries yield many
special cases. In particular, one can easily adjust to our more
general setting the examples given in Bairamov et al. (2005). We
present here only two examples making use of Corollary~2.

{\bf Example 1 (Weibull distribution).}\  Let $l_G=0$,
$r_G=\infty$, and $g(y)=y^{-\alpha(k+1)}$. Then, according to
Corollary~2, $Y$ has the Weibull distribution with $
G(y)=1-\exp\{-cy^{\alpha}\} $  if and only if for $1\le k\le n-1$
\[
E\left[[Y(n)]^{-\alpha (k+1)}\ |\ Y(n-k)=s, Y(n+1)=t\right]=
t^{-\alpha k}s^{-\alpha },
 \]
where $0<s<t<\infty$. In particular, a random variable $\tilde{Y}$
has the Inverse Weibull distribution with
$\tilde{G}(y)=\exp\{-cy^{1/2}\}$ if and only if
\[
E\left[\tilde{Y}(n)\ |\ \tilde{Y}(n-1)=s, \tilde{Y}(n+1)=t\right]=
\sqrt{st}.
 \]

 {\bf Example 2 (Pareto distribution).}\ Let
$l_G=a>0$, $r_G=\infty$, and $g(y)=[\log y]^{-(k+1)}$. Then, $Y$
has the Pareto distribution with $G(y)=1-(a/y)^c \quad (y\ge a)$
if and only if for $1\le k\le n-1$,
\[
E\left[[\log Y(n)]^{-(k+1)}\ |\ Y(n-k)=s, Y(n+1)=t\right] = [\log
t]^{-k} [\log s]^{-1},
\]
where $a\le s<t<\infty$. Note that the regression relation is
independent of $a$.

\section*{Acknowledgements}
The authors are grateful to the referees for careful reading and
constructive suggestions which were helpful in improving
substantially the presentation.

\section*{References}

\noindent Ahsanullah, M. (2004) Record values - theory and
applications. Lanham, MD:

University Press of America

\noindent Ahsanullah, M. and Raqab, M. (2006) Bounds and
Characterizations of

Record Statistics. Nova Science Publishers, New York.

 \noindent Bairamov, I., Ahsanullah, M. and Pakes, A.
(2005) A Characterization of

continuous distributions via regression on pairs of record values.
Aust. N.Z.

J. Stat. 47:543-547

\noindent Gupta, R.C. and Ahsanullah M. (2004) Some
characterization results based

on the conditional expectation of a function of non-adjacent order
statistic

(record value). Ann. Inst. Statist. Math. 721-732

\noindent Nevzorov, V.B. (2001) Records: mathematical theory.
Providence, RI:

American Mathematical Society

\noindent Pakes, A.G. (2004) Product integration and
characterization of probability

laws. J. Appl. Statist. Sci. 13:11-31

\end{document}